\documentclass[12pt,thmsa]{amsart}

\def\Cal{\mathcal}

\def\B{{\Cal B}}
\def\C{{\Cal C}}

\def\R{{\Cal R}}
\def\M{{\Cal M}}

\def\D{{\Cal D}}
\def\S{{\Cal S}}

\def\K{{\Cal K}}
\def\I{{\Cal I}}

\def\gnk{G_{n,k}}

\def\bbk{{\Bbb K}}

\def\bbr{{\Bbb R}}

\def\bbc{{\Bbb C}}

\def\cl{{\hbox{\rm cl}}}

\def\vol{{\hbox{\rm vol}}}
\def\min{{\hbox{\rm min}}}

\def\gnk{G_{n,k}}

\def\part{\partial}
\def\intl{\int\limits}
\def\b{\beta}

\def\lng{\langle}
\def\rng{\rangle}
\def\Gam{\Gamma}

\def\a{\alpha}
\def\om{\omega}

\def\del{\delta}
\def\vp{\varphi}

\def\gam{\gamma}

\def\sig{\sigma}
\def\lam{\lambda}

\newtheorem{theorem}{Theorem}[section]
\newtheorem{lemma}[theorem]{Lemma}
\newtheorem{definition}[theorem]{Definition}

\newtheorem{corollary}[theorem]{Corollary}

\theoremstyle{remark}
\newtheorem{remark}[theorem]{Remark}

\numberwithin{equation}{section}

\newcommand{\be}{\begin{equation}}
\newcommand{\ee}{\end{equation}}

\newcommand{\bea}{\begin{eqnarray}}
\newcommand{\eea}{\end{eqnarray}}
\newcommand{\Bea}{\begin{eqnarray*}}
\newcommand{\Eea}{\end{eqnarray*}}

\begin{document}

\title[Generalized cosine transforms]
{Generalized cosine transforms and classes of star bodies}

\author{Boris Rubin}
\address{
Department of Mathematics, Louisiana State University, Baton Rouge,
LA, 70803 USA\ and Einstein Institute of Mathematics, The Hebrew
University of Jerusalem, Jerusalem, 91904, Israel}

\email{borisr@math.lsu.edu}


\subjclass[2000]{Primary 44A12; Secondary 52A38}



\keywords{ Radon transforms, cosine transforms, intersection bodies}

\begin{abstract}
The spherical Radon transform on the unit sphere  in $\bbr^n$ can be
regarded  as a member of the analytic family of
  suitably normalized generalized cosine transforms.
We derive new formulas for these transforms and apply them to study
classes of intersections bodies in convex geometry. In particular,
we show that some known classes of intersection bodies are
subclasses of a more general class $\K_{\a,n}$ of origin-symmetric
star bodies in $\bbr^n$ that can be defined and characterized in
terms of the generalized cosine transforms.
\end{abstract}

\maketitle

\section{Introduction}

\setcounter{equation}{0}

This article has two sources. The first one is the theory of the
spherical Radon transforms, that amounts to classical works by H.
Minkowski, P. Funk, and S. Helgason \cite{He}. The simplest example
is the Minkowski-Funk transform, which integrates  functions  on the
unit sphere $ S^{n-1}$ in $\bbr^n$ over great circles of codimension
$1$. The spherical Radon transform can be regarded as a member of
the analytic family of
  suitably normalized generalized cosine transforms. The latter, without naming,
  first appeared in the paper  by V.I. Semyanistyi \cite{Se}
(codimension $1$) and extended by the author \cite{R1} (codimension
$\ge 1$). The analytic family of generalized cosine transforms
 (usually, without naming) has proved to be important in PDE,
 harmonic analysis, and other areas; see
\cite{Es}, \cite{Pl}, \cite{R2},  \cite{RZ}, \cite{Sa1}, \cite{Sa2},
\cite{Str1}, and references therein. Higher rank modifications of
cosine transforms were considered in \cite{Al}, \cite{AB},
\cite{GH1}, \cite{GH2}, \cite{OR}. We note that the name
``spherical" or "circular"  Radon transform is attributed in some
publications to Radon-like transforms   of different type  (see,
e.g., \cite{A}, \cite{AK},  \cite{Q}).

Another source of our article is convex geometry, related problems
in probability, stochastic geometry, and Banach space theory; see
\cite{BL}, \cite{G}, \cite{GZ}, \cite{K3}, \cite{Schn}, and
references therein. Here the name {\it cosine transform} was adopted
for the integral operator $$(\C f) (\theta)=\intl_{S^{n-1}} f(u)
|\theta \cdot u| du$$ according to pioneering works by W. Blaschke,
A.D. Alexandrov, and P. L\'evy. A more general
 {\it $p$-cosine transform}
\be\label{pcst} (\C_p f) (\theta)=\intl_{S^{n-1}} f(u) |\theta \cdot
u|^p\, du, \qquad \theta \in S^{n-1},\ee  is also commonly in use
 and reflects a number of  important geometric concepts.

In fact, both sources are intimately connected, and analytic
families associated to the spherical Radon transform include
$p$-cosine transforms up to normalization.

In the present paper, we continue our study of the generalized
cosine transforms started in \cite{R3}-\cite{R2}, keeping in mind
applications to convex geometry. Section 2 contains preliminaries.
In Section 3, we derive  new formulas, which reveal interrelation
 of different analytic families of intertwining operators on the
sphere. One of such formulas is a factorization of the
Minkowski-Funk transform as a product of mutually orthogonal
spherical Radon transforms of codimension greater than $1$; see
Theorem \ref{l3}. Section 4 deals with applications. Using results
of Section 3, we give alternative proof to some known facts in
convex geometry with the main focus on classes of intersection
bodies. We show that some of these classes, studied separately in a
series of publications, are, in fact, subclasses of a certain more
general class $\K_{\a,n}$ of origin-symmetric star bodies  that can
be characterized in terms of the generalized cosine transforms.

Our approach can also be applied
 to a series of problems related to projection bodies, $p$-centroid
 bodies, and their polars; see \cite{YY} and \cite{K3} regarding this circle of
 problems and further references.

One should also mention important works of J. Bourgain, S. Campi, R.
Gardner, P. Goodey, E. Grinberg, H. Groemer, S. Helgason, A.
Koldobsky, E. Lutwak, R. Schneider, R. Strichartz, W. Weil, G.
Zhang,  and many others, containing substantial contribution
 to harmonic analysis on the sphere in the context of
its application to integral geometry. Our list of references is far
from being complete and gives only key directions. Our interest to
this research was stimulated in part by the Busemann-Petty type
problems \cite{BZ}, \cite{K3}, \cite{RZ}, which reveal a remarkable
interplay between harmonic analysis, convex geometry, and Radon
transforms.

\section{Preliminaries}

\setcounter{equation}{0}

\subsection{Harmonic analysis on the sphere}
The main references are \cite{M}, \cite{Ne}, \cite{SW}, and a survey
article \cite{Sa2}. We use the following notation:
 $S^{n-1}$ is the unit sphere in $\bbr^n$, $C(S^{n-1})$  the spaces
 of continuous functions on $S^{n-1}$, $\, \sig_{n-1}=
2\pi^{n/2}/\Gam (n/2) $  the area of $S^{n-1}$. For $\theta \in
S^{n-1}$, $d\theta$ denotes the normalized induced Lebesgue measure
on $S^{n-1}$ and $d(\cdot, \cdot)$ stands for  the geodesic distance
on $S^{n-1}$. We denote by $ \, e_1, e_2, \ldots , e_n$  the
coordinate unit vectors; $SO(n)$ is the special orthogonal group of
$\bbr^n$; $ \, SO(n-1)$ is
 the subgroup of $SO(n)$ preserving $e_n$. For $\gam \in SO(n)$, we denote by  $d\gam $ the
 normalized $SO(n)$-invariant measure on $SO(n)$ with total mass
$1$. We use the notation
 $\D=\D(S^{n-1})$ for the space of infinitely differentiable test functions on $S^{n-1}$
equipped with the standard topology, and denote by
$\D'=\D'(S^{n-1})$ the corresponding dual space of distributions.
The subspace of even test functions (distributions) is denoted by
$\D_{e}$ ( $\D_{e}'$). The notation $\M(S^{n-1})$ is adopted for the
space of finite Borel measures on $S^{n-1}$. If $i$ is an integer,
$1\le i\le n-1$, then  $G_{n,i}$ denotes the Grassmann manifold of
$i$-dimensional linear subspaces $\xi$ of $\Bbb R^n$; $d\xi$ stands
for the normalized $SO(n)$-invariant measure on $G_{n,i}$;
$\D(G_{n,i})$ is the space of infinitely differentiable functions on
$G_{n,i}$.

Let $\{ Y_{j, k} (\theta) \}$ be an orthonormal basis of spherical
harmonics on $S^{n-1}$. Here $j = 0, 1, 2, \dots ,$ and $k = 1, 2,
\dots, d_n (j)$ where $d_n (j)$  is the dimension of the subspace of
spherical harmonics of degree $j$. Each test function $\om \in \D$
admits a decomposition $\om (\theta)=\sum_{j,k} \om_{j,k}Y_{j,
k}(\theta)$ with the Fourier-Laplace coefficients
$\om_{j,k}=\int_{S^{n-1}}\om (\theta)Y_{j, k}(\theta)d\theta$, which
decay rapidly as $j\to \infty$. Each distribution $f \in \D'$ can be
defined by $(f,\om)=\sum_{j,k} f_{j,k}\om_{j,k}$ where
 $f_{j,k}=(f,Y_{j, k})$ grow not faster than $j^m$ for some
integer $m$.

 The Poisson integral of a function $f\in L^1 (S^{n-1})$ is
defined by \be\label{pu} (\Pi_t f)(\theta)= (1-t^2) \intl_{S^{n-1}}
f(u) |\theta-t u|^{-n} du, \quad  0<t <1,\ee with the
Fourier-Laplace decomposition $\Pi_t f=\sum_{j,k} t^j f_{j,k}Y_{j,
k}$. For $f \in \D'$, this decomposition serves as a definition of
$\Pi_t f$. The space $\D (\D_e)$ is dense in $\D' (\D'_e)$ because
each distribution $f$ can be approximated in the weak sense by its
Poisson integral $\Pi_t f$, when  $t \to 1$.

A distribution $f \in \D'$ is  nonnegative if $(f,\om) \ge 0$ for
every nonnegative test function $\om$. Given a certain space $A(X)$,
consisting of functions, measures, or distributions on $X$, we
denote by $A_+(X)$ the relevant subspace of all nonnegative elements
of $A(X)$.

The following statement is a spherical analog of the well-known fact
for distributions on $\bbr^n$ \cite {Schw}. For the sake of
completeness, we present it with proof.
\begin{theorem}\label{sm}
A distribution $f \in \D'(S^{n-1})$ is nonnegative if and only if it
is a nonnegative finite measure on $S^{n-1}$, i.e., $f\in \M_+
(S^{n-1})$.
\end{theorem}
\begin{proof} The ``if" part is obvious. The proof of the `` only if"
part relies on the following

\noindent {\bf Proposition.} {\it A distribution $f \in \D'$ is  a
finite measure on $S^{n-1}$ if and only if the order of $f$ equals
$0$, i.e.,} \be\label{di1}|(f,\om)|\le c\, ||\om||_{C(S^{n-1})}
\qquad \forall \om \in \D. \ee

\noindent {\it Proof of the Proposition.} The ``only if" part is
obvious. Conversely, let (\ref{di1}) hold. Since $D$ is dense in
$C(S^{n-1})$, then  $f$ extends as a linear continuous functional on
$C(S^{n-1})$. By the Riesz theorem, there is a measure $\mu$ on
$S^{n-1}$ such that $(f,\om)=\int_{S^{n-1}}\om (\theta) d\mu
(\theta)$ for every $\om \in C(S^{n-1})$. This gives the statement.

Now we conclude the proof of the theorem. For  $\om \in \D$,
$$-||\om||_{C(S^{n-1})} \le \om\le ||\om||_{C(S^{n-1})}.$$ Hence, if
$f \in \D'$ is nonnegative, then
$$-(f,1)\,||\om||_{C(S^{n-1})} \le (f,\om)\le (f,1)\,||\om||_{C(S^{n-1})},$$
i.e., $|(f,\om)|\le (f,1)\,||\om||_{C(S^{n-1})}$ for every $\om \in
\D$. This means that $f$ has order $0$ and, by Proposition, $f$ is a
(nonnegative) finite measure.
\end{proof}

\subsection{Spherical Radon transforms}

 For continuous functions $f(\theta)$ on $S^{n-1}$ and $\varphi
(\xi)$ on $G_{n,i}$, the totally geodesic Radon transform $R_i f$
 and its
dual  $R_i^*\varphi$ are defined by \be\label{rts}
 (R_i f)(\xi) = \intl_{S^{n-1}\cap\xi} f(\theta) \, d_\xi \theta, \qquad
  (R_i^* \varphi)(\theta) = \intl_{\xi \ni \theta}  \varphi (\xi)  \, d_\theta \xi,
\ee where $d_\xi \theta$ and $ d_\theta \xi$ denote the normalized
 induced
 measures on the corresponding manifolds $S^{n-1}\cap\xi$
and $\{\xi\in G_{n,i}: \xi \ni \theta \}$; see \cite{He}, \cite{R1}.
The precise meaning of the second integral is
 \be\label{drt}
 (R_i^*
 \vp)(\theta)=\int_{SO(n-1)}\vp (r_\theta \gam p_0) \,  d\gam, \qquad \theta \in
 S^{n-1},
 \ee
 where $p_0 =\bbr
e_{n-i+1} + \ldots + \bbr e_{n}$ is the coordinate $i$-plane and $
r_\theta \in SO(n)$ is a rotation satisfying $r_\theta e_n =\theta$.
The corresponding duality relation has the form \be\label{dual}
\intl_{G_{n,i}} (R_if)(\xi) \vp (\xi) d \xi =
 \intl_{S^{n-1}} f(\theta) (R_i^* \vp) (\theta)
d\theta. \ee It is applicable when the integral in either side is
finite for $f$ and $\vp$ replaced by $|f|$ and $|\vp|$,
respectively.

The Radon transform $R_i$ and its dual extend as linear bounded
operators from $L^1(S^{n-1})$ to $L^1(G_{n,i})$ and from
$L^1(G_{n,i})$ to $L^1(S^{n-1})$, respectively. For
 finite Borel measures $\mu$ on $S^{n-1}$  and $\nu$ on $G_{n,i}$, owing to
   (\ref{dual}), we define $R_i \mu \in \M(G_{n,i})$ and $R_i^* \nu \in \M(S^{n-1})$ by
  the following equalities:
\be\label{rtm} \intl_{\gnk} \! \!(R_i \mu)(\xi) \, \vp (\xi)   d \xi
\!=  \!  \intl_{S^{n-1}} \!\! (R_i^* \vp) (\theta)   d\mu (\theta),
\quad \vp \! \in \! C(G_{n,i});\ee \be\label{drtm}
 \intl_{S^{n-1}} \!\! (R_i^* \nu) (\theta)
f(\theta) \ d\theta  \!= \!
  \intl_{G_{n,i}} \!\!
(R_i f)(\xi)   d\nu (\xi), \quad f  \!\in  \!C(S^{n-1}).\ee We also
write (\ref{dual}), (\ref{rtm}), and  (\ref{drtm}) briefly as
$$(R_if, \vp)=(f,R_i^* \vp), \quad (R_i \mu, \vp)=(\mu,R_i^* \vp),
\quad (R_i^*\nu, f)=(\nu,R_if).$$
 If $i=n-1, \; u\in S^{n-1}$,  and $\xi=u^\perp\in
G_{n, n-1}$,  it is convenient to use another notation $(R_{n-1}
f)(u^\perp)=(M f)(u)$ where \be\label{mf} (M f)(u)=\intl_{\{\theta
\,: \, \theta \cdot u =0\}} f(\theta) \,d_u\theta, \qquad u \in
S^{n-1}, \ee is the Minkowski-Funk transform of $f$.  Here
$d_u\theta$ denotes the corresponding
 normalized  measure.

\section{Analytic families}

\setcounter{equation}{0}
\subsection{Definitions and basic properties}
We start by reviewing some  facts from \cite {R1}, \cite {R2}. Given
a subspace $\xi \in G_{n,i}$, we denote by  $\text{\rm
Pr}_{\xi^\perp} \theta $
 the  orthogonal
 projection of $\theta \in S^{n-1}$ onto $\xi^\perp$, the orthogonal complement of $\xi$. Then
 $|\text{\rm Pr}_{\xi^\perp} \theta |=\sin [d(x, S^{n-1} \cap \xi)]$ is the length of
 $\text{\rm Pr}_{\xi^\perp} \theta$. We
consider analytic families of intertwining operators defined for
$f\in L^1(S^{n-1})$ and $\vp \in L^1(G_{n,i})$
 by \be\label{rka} (R_i^\a f)(\xi)=\gam_{n,i}(\a)\, \intl_{S^{n-1}}
|\text{\rm Pr}_{\xi^\perp} \theta|^{\a+i-n} \, f(\theta) \, d\theta,
\ee
 \be\label{rkad}
 (\overset *  R {}^\a \vp)(\theta)=\gam_{n,i}(\a)\,\intl_{G_{n,i}}  |\text{\rm
Pr}_{\xi^\perp} \theta|^{\a+i-n} \, \vp (\xi)\, d\xi,\ee
$$\gam_{n,i}(\a)=\frac{ \sig_{n-1}\,\Gamma((n - \a- i)/2)} {2\pi^{(n-1)/2} \,
\Gamma(\a/2)}, \qquad Re \, \a
>0, \quad \a+i-n \neq 0,2,4, \ldots .$$
For $i=n-1$, we write (\ref{rka}) as  \be\label{af} (M^\a f)(u)=
\gam_n(\a)\, \intl_{S^{n-1}} f(\theta) |\theta \cdot u|^{\a-1}
\,d\theta, \ee
\[  \gam_n(\a)={ \sig_{n-1}\,\Gamma\big(
(1-\a)/2\big)\over 2\pi^{(n-1)/2} \Gamma (\a/2)}, \qquad Re \, \a
>0, \qquad \a \neq 1,3,5, \ldots .\]
Operators (\ref{rka}) and (\ref{rkad}) were introduced in \cite{R1}
as generalizations of  (\ref{af}). The latter was introduced by
Semyanistyi \cite{Se} and studied in numerous publications; see
\cite{R2}, \cite{Sa1}, \cite{Sa2},
 and references therein. All these
operators are intimately related to the Radon transforms (\ref{rts})
and (\ref{mf}). Namely, if $f$ and $\vp$ are continuous functions,
then \cite{R1} \bea\label{lim} \lim\limits_{\a \to 0} R_i^\a
f&=&R_i^0 f =c_i \,R_if,
\qquad c_i=\frac{\sig_{i-1}}{2\pi^{(i-1)/2}};\\
\label{limd} \lim\limits_{\a \to 0} \overset *  R_i {}^\a \vp
&=&\overset *  R_i {}^0 \vp =c_i\, R_i^* \vp,\\
\label{lim1} \lim\limits_{\a \to 0} M^\a f&=&M^0 f= c_{n-1}\,
 Mf, \qquad c_{n-1}=\frac{\sig_{n-2}}{2\pi^{(n-2)/2}}.\eea
 This means that the Radon transform, its dual, and the Minkowski-Funk
transform can be regarded (up to a constant multiple) as  members of
the corresponding analytic families $\{ R_i^\a\}$, $\{\overset *
R_i {}^\a\}$, $\{ M^\a\}$.

Integrals  (\ref{rka}) - (\ref{af}) are absolutely convergent if $Re
\, \a>0$ for any integrable functions $f$ and $\vp$ . When $f$ and
$\vp$ are infinitely differentiable, these integrals extend to all
$\a \in \bbc$ as meromorphic functions of $\a$. For (\ref{af}), this
extension can be realized in terms of spherical harmonic
decomposition. Namely (see, e.g., \cite{R3}, \cite{R2}), if $f \in
\D (S^{n-1})$, then \be\label{dec} M^\a f=\sum\limits_{j, k}
m_{j,\a} f_{j, k} Y_{j, k} \ee where $f_{j, k}= \int_{S^{n-1}}
f(\theta) Y_{j, k}(\theta) d\theta$,   \be\label{mlt} m_{j, \a}
=\left\{
\begin{array}{cl} (-1)^{j/2}\, \frac{\Gamma
(j/2+(1-\a)/2)}{ \Gamma (j/2+(n-1+\a)/2)}
  &  \mbox{\rm if $j$ is even}, \\
0 &  \mbox{\rm if $j$ is odd}.
\end{array}
\right. \ee If $f\in \D'$, then, $M^\a f$ is a distribution defined
by $$(M^\a f,\om)=(f,M^\a\om)=\sum\limits_{j,k} m_{j,\a} \,f_{j,
k}\,\om_{j, k},\quad \om \in \D; \quad\a \neq 1,3,5, \ldots \,.$$

\begin{remark} The normalization in (\ref{af}) is motivated by
the asymptotic relation \be\label{as} m_{j, \a} \sim
(-1)^{j/2}(j/2)^{-(\a+n/2 -1)}, \qquad j \to \infty, \quad \text{$j$
even},\ee according to which $M^\a$ is a smoothing operator of order
$Re \,\a+n/2 -1$. It means that in many aspects, $M^\a$  acts as an
integral operator if $Re \, \a \ge 1-n/2 $ (even outside of the
domain of  absolute convergence) and as a differential operator
otherwise. Action of $M^\a$ in different scales of function spaces
(H\"older spaces, $L^p$-spaces, Sobolev spaces) was studied in
\cite{Str1}, \cite{Str2}, \cite{Sa1}, \cite[Section 2]{R4}. In
numerous publications related to integral geometry, operators
(\ref{af}) with $\a-1$ replaced by $p$ are called the (normalized)
$p$-cosine transforms (cf. (\ref{pcst})) in view of close connection
with isometric embeddings of normed spaces $(\bbr^n , || \cdot||)$
into $L^p$-spaces.
\end{remark}

The following obvious consequence of (\ref{mlt}) was widely used in
diverse publications related to the analytic family $\{ M^\a\}$; see
\cite {R2}.
\begin{lemma}\label{l1} Let $\a, \b \in \bbc; \; \a, \b \neq
1,3,5, \ldots \,$. If $\a+\b=2-n$,  then \be\label{st}M^\a
M^{\b}=I\quad \text{\rm (the identity operator)}.\ee If $\, \a,
2-n-\a \neq 1,3,5, \ldots $, then $M^\a$ is an authomorphism of the
spaces $\D_e(S^{n-1})$ and $\D'_e(S^{n-1})$.
\end{lemma}
\begin{proof} The equality $M^\a M^{\b}=I$ is
equivalent to $m_{j, \a} m_{j, \b}=1$, $\;\a+\b=2-n$. The latter
immediately follows
 from (\ref{mlt}). The second statement  is a
consequence of the standard theory of spherical harmonics \cite{Ne},
because the Fourier-Laplace multiplier $m_{j,\a}$  has a power
behavior as $j \to \infty$.
\end{proof}
\begin{corollary} The Minkowski-Funk transform $M$ on the space $\D_e(S^{n-1})$
 can be inverted by the formula
\be\label{mmm} (M)^{-1}=c_{n-1}\,M^{2-n}, \qquad
c_{n-1}=\frac{\sig_{n-2}}{2\pi^{(n-2)/2}}.\ee
\end{corollary}

Note that there is a wide variety of diverse inversion formulas for
the Minkowski-Funk transform (see \cite{He}, \cite{R2} and
references therein), but all of them
 are, in fact, different realizations of (\ref{mmm}), depending on
 classes of functions.
\begin{lemma}\label{lcon} Let $\a, \b \in \bbc; \; \a, \b \neq
1,3,5, \ldots \,$. If $Re\, \a >Re\, \b$, then $M^\a=M^\b A_{\a,\b
}$, where $A_{\a,\b }$ is a smoothing operator of order $Re
\,(\a-\b)$ with the Fourier-Laplace multiplier \be\label{ab}
a_{\a,\b }(j)= \frac{\Gamma (j/2+(1-\a)/2)}{ \Gamma
(j/2+(n-1+\a)/2)}\, \frac{\Gamma (j/2+(n-1+\b)/2)}{ \Gamma
(j/2+(1-\b)/2)},\ee so that $a_{\a,\b }(j)\sim (j/2)^{\b-\a}$ as $ j
\to \infty$. If $\a$ and $ \b$ are real numbers satisfying $\a>
\b>1-n,\; \a+\b<2$, then $A_{\a,\b }$ is
 an integral operator with the property $A_{\a,\b }f\ge 0$ for any
 nonnegative  $f \in L^1 (S^{n-1})$.
\end{lemma}
\begin{proof} The first statement follows from (\ref{mlt}). To prove
the second one, we consider  integral operators \bea\label{q1}\quad
(Q_{+}^{\mu,\nu} f)(x) &=& \frac{2}{\Gam (\mu/2)} \int_0^1
(1 - t^2)^{\mu/2 -1} (\Pi_t f)(x) \, t^{n-\nu} dt,\\
\label{q2}\quad (Q_{-}^{\mu,\nu} f)(x) &=& \frac{2}{\Gam (\mu/2)}
\int_1^\infty (t^2-1)^{\mu/2 -1} (\Pi_{1/t} f)(x) \, t^{1-\nu}
dt,\eea containing the Poisson integral (\ref{pu}). The
Fourier-Laplace multipliers of $Q_{+}^{\mu,\lambda}$ and
$Q_{-}^{\mu,\nu}$ are
 \be\label{mqa} \hat q_{+}^{\mu,\nu}(j)\!=
\!\frac{\Gam((j\!+\!n\!-\!\nu \!+\!1)/2)}{\Gam((j\!+\!n\!-\!\nu
\!+\!1\!+\!\mu)/2)}, \quad
 \hat  q_{-}^{\mu,\nu}(j)\!=\!
\frac{\Gam((j\!+\!\nu\! -\!\mu)/2)}{\Gam((j\!+\!\nu )/2)}.\ee They
can be easily calculated by taking into account that $\Pi_t \sim
t^j$ in the Fourier-Laplace terms \cite[ p. 145]{SW}. Integrals
(\ref{q1}) and (\ref{q2}) are absolutely convergent for $f\in
L^1(S^{n-1})$ provided $0<Re\, \mu<Re\, \nu<n$ and preserve
positivity of $f$. Comparing (\ref{mqa}) and (\ref{ab}), we obtain
a factorization  $A_{\a,\b }=Q_{+}^{\a-\b,1-\b} Q_{-}^{\a-\b,1-\b}$
(set $\mu=\a-\b,\; \nu=1-\b$), which implies the second statement of
the lemma.
\end{proof}

The next analytic family we deal with is the family of the
generalized sine transforms \bea\label{sin} \qquad (Q^\a
f)(\theta)\!&\!\!=\!\!&\!\frac{\sig_{n-1}\Gam ((n\!-\!1\!-\!\a)/2) }
{2 \pi^{(n-1)/2} \Gam (\a/2)}\!\intl_{S^{n-1}}\!\!(1\!-\!|u \cdot
\theta|^2)^{(\a-n+1)/2} f(u) du
\\\!&\!\!=\!\!&\!\frac{\sig_{n-1}\Gam ((n\!-\!1\!-\!\a)/2) } {2
\pi^{(n-1)/2} \Gam (\a/2)}\! \intl_{S^{n-1}}\!\!(\sin [d(u,
\theta)])^{\a-n+1} f(u) du,\nonumber \eea
 where $Re \,\a > 0, \quad \a-n \neq 0,2,4, \ldots \, $. Operators $Q^\a$
 serve as  analogues of Riesz potentials in the theory of
spherical Radon transforms  \cite{He}, \cite{Str2}. Detailed
investigation of operators (\ref{sin}), including inversion
formulas, can be found in \cite{R1}.    The Fourier-Laplace
multiplier
 of $Q^\a$  has the form \be\label{muq} \hat
q_\a(j)=\frac{\Gamma \big( \frac{j+n-1-\a}{2} \big) \, \Gamma \big(
\frac{j+1}{2}\big)}{\Gamma \big( \frac{j+\a+1}{2}\big) \,
\Gamma\big( \frac{j+n-1}{2}\big)}
 \qquad ( \sim (j/2)^{-\a},  \quad j \to \infty)
 \ee
for $j$ even, and $\hat q_\a(j)=0$ for $j$ odd . For $Re \,\a \le
0$,
 $Q^\a f$ is defined by analytic continuation.  The following equalities
 play an important role in the theory of Radon transforms on $S^{n-1}$. If $f \in
 \D_{e}(S^{n-1})$  and $ \a \in \bbc,  \quad \a+i-n \neq 1,3,5, \ldots  $,
 then \cite{R1}
\be\label{rq} {\overset *  R}{}_i^\a R_i f=\lam_1 Q^{\a+i-1}f,
\qquad R^*_i R_i^\a f=\lam_2 Q^{\a+i-1}f,\ee
\[ \lam_1=\frac{\Gam((n-1)/2)}{\sig_{n-1}
\, \Gam((n-i)/2)}, \qquad \lam_2 =\frac{
\Gam((n-1)/2)}{\Gam((n-i)/2)}.\] In particular, by (\ref{lim}),
\be\label{hru1} R^*_i R_i f\!=\!c \,  Q^{i-1}f,\quad
c\!=\!\frac{2\pi^{(i-1)/2}\, \Gam((n\!-\!1)/2)}{ \sig_{i-1}\,
\Gam((n\!-\!i)/2)}\quad \left (\!=\!\frac{1}{\hat q_{i-1}(0)}\right
),\ee and \be\label{hru} M^\a M^0 f=Q^{\a+n-2}f.\ee The latter
follows directly from (\ref{muq}) and (\ref{mlt}). Furthermore, by
(\ref{muq}), $ \, Q^0 f =f \, $, and (\ref{rq}) yields the following
inversion formula: \be\label{inr} {\overset * R}{}_i^{1-i} R_i
f=\lam_1 f. \ee

The next statements  represent main results of the ``analytic part"
of the paper. We establish new connections between operator families
defined above. For $\xi \in G_{n,i}$, we denote \be\label{prp}
(R_{n-i,\perp}f)(\xi)=(R_{n-i}f)(\xi^{\perp}), \qquad
(R^\a_{n-i,\perp}f)(\xi)=(R^\a_{n-i}f)(\xi^{\perp}).\ee
\begin{lemma}\label{l2} Let $f \in L^1(S^{n-1}), \quad Re\, \a >0; \quad \a \neq 1,3,5, \ldots \,
$. Then \be \label{con} (R_i M^\a f)(\xi)=
\frac{2\pi^{(i-1)/2}}{\sig_{i-1}} \, (R_{n-i,\perp}^{\a +i-1}
f)(\xi), \qquad \xi \in G_{n,i},\ee or (replace $i$ by
 $n-i$)\be \label{conn} (R_{n-i,\perp} M^\a f)(\xi) = \frac{2\pi^{(n-i-1)/2}}{\sig_{n-i-1}} \,
(R_i^{\a +n-i-1} f)(\xi).\ee If
 $f \in \D_e(S^{n-1})$, then (\ref{con}) and (\ref{conn}) extend to $Re\, \a \le 0$ by
 analytic continuation. In particular,
  \be\label{kr} (R_i M^{1-i}
  f)(\xi)=\tilde c \,(R_{n-i, \perp}f)(\xi), \quad \tilde c=
\frac{\sig_{n-i-1}\,\pi^{i-n/2}}{\sig_{i-1}}.\ee
\end{lemma}
\begin{proof} For $Re\, \a >0$,
\[(R_i M^\a f)(\xi)=\gam_n(\a) \int_{S^{n-1}\cap\xi} \, d_{\xi}
u  \intl_{S^{n-1}} f(\theta) |\theta\cdot u|^{\a-1} \,d\theta.\]
Since $|\theta\cdot u|=|\text{\rm Pr}_{\xi} \theta||v_\theta \cdot
u|$ for some $v_\theta \in S^{n-1}\cap\xi$, by changing the order of
integration, we obtain \[ (R_i M^\a f)(\xi)=\gam_n(\a)\,
\intl_{S^{n-1}} f(\theta) |\text{\rm Pr}_{\xi} \theta|^{\a-1}
\,d\theta  \int_{S^{n-1}\cap\xi} |v_\theta \cdot u|^{\a-1} d_{\xi}
u.\] The inner integral is independent on $v_\theta$ and can be
easily evaluated: \bea \int_{S^{n-1}\cap\xi} |v_\theta \cdot
u|^{\a-1} d_{\xi} u&=&\frac{\sig_{i-2}}{\sig_{i-1}}\intl_{-1}^1
|t|^{\a-1} (1-t^2)^{(i-3)/2}\, dt \nonumber\\
&=&\frac{2\pi^{(i-1)/2} \,\Gamma (\a/2)}{ \sig_{i-1}\,\Gamma(
(i+\a-1)/2)}.\nonumber\eea
 This implies (\ref{con}). Formula (\ref{kr}) then follows by analytic continuation (just set
$\a=1-i$ and make use of (\ref{lim}) with $i$ replaced by $n-i$).
 \end{proof}

Equality (\ref{kr}), written in terms of  Fourier transforms, was
actually  obtained by Koldobsky; see \cite[Lemma 7]{K2} and
\cite[Corollary 1]{K4}. The  argument in these works essentially
differs from ours.
\begin{lemma}\label{l33} For all $\a \in \bbc$ such that $\, \a,
2-n-\a \neq 1,3,5, \ldots $, \be\label{coi} R_i^\a (\D_e)=R_i
(\D_e).\ee Specifically, $R_i^\a f=R_i f_1$ for functions $f$ and $
f_1$  in $\D_e$ connected by \bea\label{coi1}
f&=&\frac{2\pi^{(i-1)/2}}{\sig_{i-1}}\,M^{1-n+i}M^{1-\a-i}f_1, \\
\label{coi11}f_1&=&\frac{\pi^{(1-i)/2}\, \sig_{i-1}}{2}\,
M^{1-i}M^{\a+i +1-n}f.\eea
\end{lemma}
\begin{proof} We use (\ref{conn}) with  $\a$ replaced by $\a+1+i-n$ and then apply
(\ref{kr}). This gives \bea R_i^\a f&=&\frac{\pi^{(i+1-n)/2}\,
\sig_{n-i-1}}{2}\,R_{n-i,\perp} M^{\a+i
+1-n}f\nonumber\\
&=&\frac{\pi^{(1-i)/2}\, \sig_{i-1}}{2}\, R_i M^{1-i}M^{\a+i
+1-n}f\!=\!R_i f_1.\nonumber\eea
 Equality (\ref{coi1}) then
follows from (\ref{st}).
\end{proof}

 The next statement contains an intriguing factorization of the
 Minkowski-Funk transform in terms of Radon transforms associated to mutually orthogonal subspaces.
\begin{theorem}\label{l3} For $f \in L^1(S^{n-1})$ and $0<i<n$,
\be\label{svv} Mf=R_i^* R_{n-i,\perp} f.\ee
\end{theorem}
\begin{proof} By (\ref{drt}),
\bea (R_i^* R_{n-i,\perp} f)(\theta)&=&\intl_{SO(n-1)}
(R_{n-i,\perp} f)(r_\theta \gam \bbr^i)\, d\gam\nonumber \\
&=&\intl_{SO(n-1)} (R_{n-i}
f)(r_\theta \gam \bbr^{n-i})\, d\gam \nonumber \\
&=& \intl_{SO(n-1)}d\gam \intl_{S^{n-1}\cap r_\theta \gam
\bbr^{n-i}} f(v) \, dv \nonumber \\
&=& \intl_{S^{n-1}\cap \bbr^{n-i}} dw \intl_{SO(n-1)}f(r_\theta \gam
w) \,d\gam.\nonumber \eea The inner integral is independent on $w
\in S^{n-1}\cap \bbr^{n-i}$ and equals $(M f)(\theta)$. This gives
(\ref{svv}).
\end{proof}
 The following statement reveals a remarkable interplay between different analytic
families and gives a series of explicit representations of the right
inverse of the dual Radon transform  $R_i^*$ (note that $R_i^*$ is
non-injective on $\D (G_{n,i})$ when $1<i<n-1$).
 \begin{lemma}\label{l34} For  $0<i<n$, every  function $f\in \D_e (S^{n-1})$ is
represented by the dual Radon transform  $f=R_i^* g, \quad g=Af$,
where the operator $A:\D_e (S^{n-1}) \to \D(G_{n,i})$ has the
following forms: \bea\label{crl} \qquad
(Af)(\xi)&=&\label{exx1}\frac{\sig_{n-2}}{2\pi^{n/2
-1}}\,(R_{n-i,\perp} M^{2-n}f)(\xi)\\
&=&\label{exx2}\frac{\pi^{(1-i)/2}\sig_{n-2}}{\sig_{n-i-1}}\,(R_i^{1-i}f)(\xi)\\
&=&\label{exx3}\frac{\pi^{1-i}\sig_{n-2}\,\sig_{i-1}}{2\,\sig_{n-i-1}}\,(R_i
(Q^{i-1})^{-1}f)(\xi), \qquad \xi \in G_{n,i}.\eea
\end{lemma}
\begin{proof} We first show that expressions (\ref{exx1})-(\ref{exx3})
coincide. The coincidence of (\ref{exx1}) and (\ref{exx2}) follows
from (\ref{conn}). To show that  (\ref{exx2}) coincides with
(\ref{exx3}), we set $f=Q^{i-1}f_1$. Since $Q^{i-1}$ is injective,
it suffices to check the equality
$2\pi^{(i-1)/2}R_i^{1-i}Q^{i-1}f_1=\sig_{i-1}\,R_i f_1$. The latter
holds by Lemma \ref{l33}. Indeed, for $\a=1-i$, equalities
(\ref{coi1}) and (\ref{hru}) yield
$$\sig_{i-1}\,R_i
f_1=2\pi^{(i-1)/2}R_i^{1-i}M^{1-n+i}M^{0}f_1=2\pi^{(i-1)/2}R_i^{1-i}
Q^{i-1}f_1.$$  It remains to note that the representation  $f=R_i^*
g$ with $g=Af$ defined by (\ref{exx2}) follows from the second
equality in (\ref{rq}), if we apply it with $\a=1-i$ and take into
account that $Q^0 f =f$.
\end{proof}
\begin{remark} As it was mentioned above,  the map
$R_i^*: \D(G_{n,i}) \to \D_e (S^{n-1})$ is non-injective. In fact,
every function $\vp \in \D(G_{n,i})$ is represented as a sum
$\vp=\vp_R + \vp_0$ where $\vp_R $ belongs to the range
$\R_i(\D_e(S^{n-1}))$ and $\vp_0 \in \ker\R^*_i$. Indeed, let $\vp_R
=c^{-1}\,R_i[(Q^{i-1})^{-1}R_i^* \vp]$, $c$ being a constant from
(\ref{hru1}), $ \vp_0=\vp-\vp_R$. Then, by (\ref{hru1}), \bea R_i^*
\vp_0&=&R_i^* \vp-c^{-1}\,R_i^*\vp_R =R_i^* \vp-c^{-1}\,R_i^*
R_i[(Q^{i-1})^{-1}R_i^* \vp]\nonumber
\\&=&R_i^* \vp-Q^{i-1}(Q^{i-1})^{-1}R_i^* \vp=R_i^*
\vp-R_i^*\vp=0,\nonumber
 \eea i.e., $\vp_0 \in \ker R^*_i$.
\end{remark}

The following  statement is dual to Lemma \ref{l2}.
\begin{lemma}\label{l4} Let $\vp \in L^1(G_{n,i}), \;
\vp^\perp (\eta)=\vp(\eta^\perp), \quad \eta \in G_{n,n-i}$. If
$Re\, \a
>0, \quad \a \neq 1,3,5, \ldots \, $, then
\be \label{cnd} M^\a R_i^* \vp= c\, \overset * R{}^{\a +i-1}_{n-i}
\vp^\perp, \qquad c= \frac{2\pi^{(i-1)/2}}{\sig_{i-1}}.\ee If $\vp
\in \D(G_{n,i})$, then (\ref{cnd}) extends to all complex $\a\neq
1,3,5, \ldots $ by analytic continuation. In particular, by
(\ref{limd}), \be\label{ccrl} M^{1-i} R_i^* \vp= \tilde c\,
R^*_{n-i} \vp^\perp, \qquad \tilde c= \frac{\pi^{i-n/2}\,
\sig_{n-i-1}}{\sig_{i-1}}.\ee
\end{lemma}
\begin{proof} Let $Re\, \a
>0$. Owing to (\ref{dual}) and (\ref{con}), for any
 test function $\om \in \D_e (S^{n-1})$ we have \bea
&&\intl_{S^{n-1}} (M^\a R_i^* \vp)(\theta)\, \om(\theta)\,
d\theta=\intl_{S^{n-1}}(R_i^* \vp)(\theta)\, (M^\a\om)(\theta)\,
d\theta\nonumber
\\&=&\intl_{G_{n,i}} \vp(\xi)\,
(R_i M^\a\om)(\xi)\,
d\xi=\frac{2\pi^{(i-1)/2}}{\sig_{i-1}}\,\intl_{G_{n,i}}\vp(\xi)\,
(R_{n-i,\perp}^{\a +i-1} \om)(\xi)\, d\xi\nonumber \\
&=&\frac{2\pi^{(i-1)/2}}{\sig_{i-1}}\, \intl_{G_{n,n-i}}
\vp^\perp(\eta)\,(R_{n-i}^{\a +i-1} \om)(\eta)\, d\eta\nonumber
\\&=&\frac{2\pi^{(i-1)/2}}{\sig_{i-1}}\, \intl_{S^{n-1}} (\overset
* R{}^{\a +i-1}_{n-i} \vp^\perp)(\theta)\, \om(\theta)\,
d\theta.\nonumber \eea This gives the desired result.
\end{proof}
\begin{corollary}\label{crr}  If $\vp \in L^1(G_{n,i})$,
 then for any $Re\, \a
>0$, $\; \a \neq 1,3,5, \ldots \, $,
\be\label{iff1}f=R_i^* \vp\quad \text{\rm if and only if}\quad M^\a
f=c\,\overset * R{}^{\a +i-1}_{n-i} \vp^\perp,\ee $c$ being a
constant from (\ref{cnd}).
 If $\vp \in \D(G_{n,i})$, then (\ref{iff1}) extends to all
complex $\a\neq 1,3,5, \ldots $ by analytic continuation. In
particular, for $\a=1-i$, \be\label{iff}f=R_i^* \vp\quad \text{\rm
if and only if}\quad M^{1-i} f=c\,R^*_{n-i} \vp^\perp.\ee
\end{corollary}
\begin{proof} The statement is a consequence of (\ref{cnd}) and injectivity of
$M^\a$.
\end{proof}

\section{Positive definite homogeneous distributions and star bodies}

In this section we invoke the Fourier transform of homogeneous
distributions and  apply some results of Section 3 to study classes
of star bodies arising  in convex geometry.

 \subsection{The  Fourier transform of homogeneous distributions}
This is one of the oldest topics in the theory of distributions, and
there is a vast literature on this subject; see, e.g.,  \cite{GS},
\cite{Se}, \cite {Le}.  Let $\S(\bbr^n)$ be the Schwartz space of
rapidly decreasing infinitely differentiable functions on $\bbr^n$
and $\S'=\S'(\bbr^n)$ the corresponding space of tempered
distributions. The  Fourier transform of $F \in \S'$ is  defined by
$$ \lng\hat F, \hat \phi\rng= (2\pi)^{n}\lng F, \phi\rng, \quad \hat \phi(y)=
\int_{\bbr^{n}} \phi(x) \, e^{-ix \cdot y} \, dx, \quad \phi \in
\S(\bbr^n).$$  A distribution $F\in \S'$ is
 positive definite if  the Fourier transform $\hat F$  is a
  positive distribution, i.e.,  $\lng\hat F, \phi\rng \ge 0$ for every nonnegative test
function $\phi$. Given a distribution $F \in \S'$ and a complex
number $\lam$, we say that $F$ is a homogeneous distribution of
degree $\lam$  if for any $\phi \in \S(\bbr^n)$ and any $a >0$,
$\lng F, \phi (x/a)\rng =a^{\lam +n}\,\lng F, \phi\rng$. Homogeneous
distributions on $\bbr^n$ are intimately related to distributions on
the unit sphere.  We define  a homogeneous continuation of a
function $\vp$ on $S^{n-1}$ by \be \label{el}(E_\lam
\vp)(x)=|x|^\lam \vp (x/|x|), \qquad x \in \bbr^n \setminus
\{0\}.\ee If $\lam\neq -n, -n-1, -n-2, \ldots $, then the operator
$E_\lam $  extends to distributions $f\in\D'$ by the formula
\be\label{el} \lng E_\lam f, \phi \rng =(f,\phi_\lam),\ee where
$\phi \in \S(\bbr^n), \quad \phi_\lam (\theta)=\int_0^\infty r^{\lam
+n-1} \phi(r \theta) \, dr \in \D(S^{n-1})$. For $Re\, \lam \le -n$,
$\; \lam\neq -n, -n-1, -n-2, \ldots $, the last integral is
understood in the sense of analytic continuation.
\begin{lemma}\label{ds} Let $\lam \in \bbc; \;\lam\neq -n, -n-1, -n-2, \ldots
$. Then $E_\lam $ is a linear continuous operator from $\D'$ to
$\S'$.
\end{lemma}
\begin{proof} If $f_j \in \D'$ and
$\lim\limits_{j \to \infty}(f_j, \om)=(f,\om)$ for any $\om \in \D$,
then for any $\phi \in \S(\bbr^n)$ we have $\lim\limits_{j \to
\infty} \lng E_\lam f_j, \phi \rng= \lim\limits_{j \to \infty}
(f_j,\phi_\lam)=(f,\phi_\lam)=\lng E_\lam f, \phi \rng.$
\end{proof}

 The following theorem due to Lemoine \cite{Le}
characterizes the structure of homogeneous distributions.
\begin{theorem} Let  $\tau$ be a
homogeneous distribution of degree $\lam \in \bbc$.

{\rm a)} If $\lam$ is not an integer $\le -n$, there is an $f \in
\D'$ such that  $\tau=E_\lam f.$

{\rm b)} If $\lam=-n, -n-1, \ldots$, there are $f \in \D'$ and a
polynomial $P_{-\lam -n}$ homogeneous of degree $-\lam -n$ such that
$\tau=E_\lam f+P_{-\lam -n} (D)\del$, where $D=(\partial/\partial
x_1, \ldots, \partial/\partial x_n)$ and $\del$ is the Dirac
measure.
\end{theorem}

 The operator family $\{M^\a \}$ arises in the Fourier
analysis of homogeneous distributions in a natural way thanks to the
formula \be \label{cf} [E_{1-n-\a} f]^\wedge=2^{1-\a}
\pi^{n/2}\,E_{\a-1}M^\a f.\ee This holds for a $C^\infty$ even
function $f$
 and any complex $\a$ satisfying \be\label{alf}\a \notin \{1,2,3, \ldots \}
 \cup \{1-n, -n, -n-1, \ldots \}.\ee Formula (\ref{cf})  is understood in the
$S'$-sense. Namely,
 \be\label{sv} \lng E_{1-n-\a} f, \phi\rng=
2^{1-\a} \pi^{n/2}\, \lng E_{\a-1}M^\a f,
 \hat\phi\rng, \qquad \phi \in \S(\bbr^n),\ee  where both sides are
interpreted in the sense of analytic continuation. This formula (and
 a more general one for arbitrary, not necessarily even, functions) is
 known in analysis for many years and has many applications; see, e.g., \cite{Se},
\cite{Es}, \cite{Pl}, \cite{Sa1}, \cite{Sa2}, \cite{R3}, \cite{K3},
and references therein. Generalizations of (\ref{cf}) to functions
on Stiefel and Grassmann manifolds were obtained in \cite{OR}. Since
$\D_e$ is dense in $\D'_e$, then, owing to Lemmas \ref{l1} and
\ref{ds}, equalities (\ref{cf}) and (\ref {sv}) extend to all even
distributions on $S^{n-1}$.

The following lemma makes a bridge between positive definite
homogeneous distributions and operators $M^\a$.
\begin{lemma}\label{pr0}
Let $\a$ be a real number satisfying (\ref{alf}). If $f \in \D_e$,
then $E_{1-n-\a}f$ is a positive definite distribution if and only
if $(M^\a f)(\theta) \ge 0$ for every $\theta \in S^{n-1}$. If $f
\in \D'_e$, then $ E_{1-n-\a}f$ is a positive definite distribution
if and only if $M^\a f$ is a nonnegative finite measure on
$S^{n-1}$.
\end{lemma}
\begin{proof} The ``if" part is obvious from (\ref{sv}).
Indeed, if $\mu=M^\a f$ and $\phi \in \S_+(\bbr^n)$, then, by
(\ref{el}), $\lng E_{\a-1}\mu, \phi \rng=(f,\phi_\lam)\ge 0$, and by
(\ref{cf}), $ E_{1-n-\a}f$ is a positive definite distribution.
Conversely, if $ E_{1-n-\a}f$ is positive definite, then $\lng
E_{\a-1}M^\a f, \phi \rng$ is  nonnegative for every  $\phi \in
\S_+(\bbr^n)$. Choose $\phi$ of the form $\phi (x)=\psi (|x|) \om
(x/|x|)$, where $\psi$ is a smooth positive function satisfying
$\int_0^\infty r^{\a +n-2} \psi (r) dr=1$, and $\om$ is a
nonnegative test function on  $S^{n-1}$. Then, by (\ref{el}), $\lng
E_{\a-1}M^\a f, \phi \rng =(M^\a f, \om)\ge 0$. Hence, by
 Theorem \ref{sm}, $M^\a f$ is a finite measure. To get a pointwise
 inequality $(M^\a f)(\theta) \ge 0$ in the first statement of the
 lemma, we choose $\om  (\cdot)$ to be  the Poisson kernel $p_{t,\theta} (\cdot)$
 defined by
$$ p_{t,\theta}
(u)= \frac{1-t^2}{ (1-2t u\cdot \theta +t^2)^{n/2}},\qquad  0<t <1;
\quad u, \theta \in S^{n-1}.$$ Then $(M^\a f, \om)$ is, actually,
the Poisson integral of $M^\a f$ that tends to $M^\a f$ as $t \to 1$
uniformly on $S^{n-1}$ \cite{SW}. This gives the result.
\end{proof}

 \subsection{Classes of star bodies}

 Below we introduce  classes of origin-symmetric star bodies in
$\bbr^n$ associated with the analytic family $M^\a$ of the
 generalized cosine transforms. These classes include intersection bodies and
some other
 classes of bodies commonly used in  convex geometry. Let $K$
be a compact subset of  $\bbr^n, \; n \ge 2$, star-shaped with
respect to the origin. For  $x \in \bbr^n \setminus \{0\}$, the {\it
radial function} of $K$ is defined by $\rho_K(x)=\sup \{ \lam \ge 0:
\lam x \in K\}$. If $\theta \in S^{n-1}$, then $\rho_K(\theta)$ is
the Euclidean distance from the origin to the boundary of $K$ in the
direction of $\theta$. If $\rho_K$ is a positive continuous function
on $S^{n-1}$, then $K$ is said to be a {\it star body}. We denote be
$\bbk^n$ the class of all origin-symmetric star bodies  in $\bbr^n$.
The {\it Minkowski functional}  of $K \in \bbk^n$ is defined by $
||x||_K =\min \{a \ge 0 \, : \, x \in aK\}$ so that
$||\theta||_K=\rho_K^{-1}(\theta)$. A body $K\in \bbk^n$ is called
infinitely smooth if $\rho_K\in \D_e(S^{n-1})$. We say that a
sequence of bodies $K_j\in \bbk^n$ converges to $K \in \bbk^n$ in
the radial metric if $\lim\limits_{j \to \infty}|| \rho_{K_j} -
\rho_K||_{C(S^{n-1})}=0$. Given a subset $\K \subset \bbk^n$, we
denote by $\cl \,\K$ the closure of $\K$ in the radial metric.
\begin{definition}\label{df}
Let $\a$ be a real number, \be\label{alf1} \a\notin \{0,-2,-4,
\ldots\} \cup \{n, n+2, n+4, \ldots\}.\ee We define the following
classes of origin-symmetric star bodies: \bea
\label{b1}&{}&\K_{\a,n}\!=\!\{K \in \bbk^n : \, \rho_K^\a =M^{1-\a}
\mu \; \;\text{for some} \; \;\mu
\in \M_{e+}(S^{n-1})\}; \\
 \label{b2}&{}&\K\B_{\a,n}\!=\!\{K \in
\bbk^n : \, \rho_K^\a =M^{1-\a}  \rho_L^{n-\a}\;\; \text{for some
body}\;
\; L\in \bbk^n \};\\
\label{b3}&{}&\K\B_{\a,n}^\infty\!=\!\{K \in \K\B_{\a,n}:
\,\rho_K\in \D_e(S^{n-1})\}.\eea \end{definition}

Some comments are in order:

1. The equality $\rho_K^\a =M^{1-\a} \mu$ in (\ref{b1}) means that
\be (\rho_K^\a, \om) =(\mu, M^{1-\a}\om) \qquad \forall \om \in
\D(S^{n-1}). \ee

2. If  $K\in \K\B_{\a,n}^\infty$, then $\rho_K^\a =M^{1-\a}
\rho_L^{n-\a}$ where, by (\ref{st}),  $
\rho_L^{n-\a}=M^{1-n+\a}\rho_K^\a \in \D_{e+}(S^{n-1})$.

3. For some $\a$, the class $\K_{\a,n}$ looks somewhat artificial
and does not contain  such a nice body as the unit ball $B$. Indeed,
since $\rho_B=1$, then,  owing to (\ref{mlt}),
$M^{1-n+\a}\rho_B^\a=\Gam ((n-\a)/2)/\Gam (\a/2)$. This is negative
if \be\label{apl} \a \in \left (\bigcup\limits _{k=0}^\infty (-4k-2,
-4k) \right )\bigcup\limits \left (\bigcup\limits_{k=0}^\infty
(n+4k, n+4k+2) \right ).\ee  Thus $B\notin \K_{\a,n}$ for all such
$\a$.
\begin{theorem}\label{tcl} Let $\a$ be a real number satisfying (\ref{alf1}).
 Then \be\label{cl2}
\K_{\a,n}=\cl \, \K\B_{\a,n}=\cl \,  \K\B_{\a,n}^\infty.\ee
\end{theorem}
\begin{proof} STEP 1. We first prove that $\K_{\a,n}\subset \cl \,
\K\B_{\a,n}^\infty$. Let $K \in \K_{\a,n}$, i.e., $\rho_K^\a
=M^{1-\a} \mu, \;  \mu \in \M_{e+}(S^{n-1})$. Our aim is to define a
sequence $K_j \in \K\B_{\a,n}^\infty$ such that $\rho_{K_j} \to
\rho_{K}$ in the $C$-norm. Consider the Poisson integral $\Pi_t
\rho_K^\a$ that converges to $\rho_K^\a$ in the $C$-norm when $t\to
1$. For any test function $\om \in \D$, we have $$ (\Pi_t
\rho_K^\a,\om)=(\rho_K^\a,\Pi_t \om)=(\mu, M^{1-\a}\Pi_t
\om)=(M^{1-\a}\Pi_t\mu, \om),$$ where $\Pi_t\mu \in
 \D_{e+}$. Hence, one can choose  $K_j \in \K\B_{\a,n}^\infty$ so that
 $\rho_{K_j}^\a=\Pi_{t_j}\rho_K^\a=    M^{1-\a}\Pi_{t_j}\mu$, where $t_j$ is a sequence in
 $(0, 1)$ approaching $1$. Clearly, $K_j $  converges to
 $K$ in the radial metric.

Conversely, let $K \in \cl \,  \K\B_{\a,n}^\infty$. It means that
there is a sequence of bodies $K_j \in \K\B_{\a,n}^\infty$ such that
$\lim\limits_{j \to \infty} ||\rho_K - \rho_{K_j}||_{C}=0$ and
$\rho_{K_j}^\a=M^{1-\a}\rho_{L_j}^{n-\a}$, $ \rho_{L_j}\in \D_{e+}$.
Then $\rho_{K_j}^\a$ approaches $\rho_{K}^\a$ in the $C$-norm, and
for any $\om\in\D_{e+}$, the expression $(\rho_{K_j}^\a,
M^{1-n+\a}\om)$ is nonnegative because
$$(\rho_{K_j}^\a,
M^{1-n+\a}\om)=(M^{1-\a}\rho_{L_j}^{n-\a},M^{1-n+\a}\om)=
(\rho_{L_j}^{n-\a},\om)\ge 0.$$ If $j \to \infty$, then
$$(\rho_{K_j}^\a,
M^{1-n+\a}\om) \to (\rho_{K}^\a, M^{1-n+\a}\om)=(
M^{1-n+\a}\rho_{K}^\a, \om) \ge 0.$$ Hence, by Theorem \ref{sm},
$M^{1-n+\a}\rho_{K}^\a$ is a nonnegative measure. Let
$\mu=M^{1-n+\a}\rho_{K}^\a$. By (\ref{st}),  for any $\om \in \D$ we
have
$$ (\rho_{K}^\a,
\om)=(M^{1-n+\a}\rho_{K}^\a, M^{1-\a} \om)=(\mu,M^{1-\a}
\om)=(M^{1-\a}\mu, \om),$$ i.e., $K \in \K_{\a,n}$. This gives
$\K_{\a,n}=\cl \,  \K\B_{\a,n}^\infty$.

STEP 2. Let us prove that $\K_{\a,n}=\cl \,  \K\B_{\a,n}$. Since
$\K\B_{\a,n}^\infty \subset\K\B_{\a,n}$, then, by Step 1,
$\K_{\a,n}\subset \cl \,  \K\B_{\a,n}$. The proof of the reverse
inclusion coincides with  the second part in Step 1 with the only
difference that now $\rho_{L_j}$ are only continuous and not
necessarily smooth.
\end{proof}
\begin{theorem}\label{bpd} A  body $K\in \bbk^n$ belongs to  $\K_{\a,n}$
  if and only if
 $||\cdot||_K^{-\a}$ is a positive definite distribution.
\end{theorem}
\begin{proof}
By Lemma \ref {pr0}, $||\cdot||_K^{-\a}\equiv E_{-\a} \rho_K^\a$ is
a positive definite distribution if and only if there is a
 measure $\mu\in \M_{e+}(S^{n-1})$ such that
$(M^{1-n+\a}\rho_K^\a, \om)=(\mu,\om)$ $\; \forall \om \in \D_e$.
Owing to (\ref{st}), the latter is equivalent to $(\rho_K^\a,
 \tilde \om)=(\mu, M^{1-\a}\tilde \om) \; \forall \tilde \om \in \D_e$
  (choose $\om=M^{1-\a}\tilde \om$). This is what we need.
\end{proof}

 Below we consider some examples when known geometric objects are members
  of the class $\K_{\a,n}$.   If $K\in
\bbk^n$ and $\xi$ is an $i$-dimensional subspace of $\bbr^n$, i.e.,
$\xi \in G_{n,i}, \; 1\le i< n$, then the volume of the
cross-section $K\cap \xi$ can be evaluated by \be\label{voli}
\vol_i(K\cap \xi) = \frac{\sig_{i-1}}{i}\intl_{S^{n-1}\cap \xi}
\rho_K^i (\theta) \,d_\xi \theta=\frac{\sig_{i-1}}{i}\,(R_i \rho_K^i
)(\xi) .\ee This can be easily obtained by passing to polar
coordinates in $\xi$.

{\bf Example 1.}\label{ex1} According to Lutwak \cite{Lu}, the class
$\I\B_n$ of {\it intersection bodies of star bodies} $\,$ in
$\bbr^n$ is defined as the range of the map $\I\B: \bbk^n \to
\bbk^n$ by the rule
$$\rho_{\I\B(L)} (\theta)=\vol_{n-1}(L\cap \theta^\perp), \qquad
\theta \in S^{n-1}, \quad L \in \bbk^n.$$ By (\ref{voli}),  it means
that $K=\I\B(L)$ if and only if \be\label{eqv2}
\rho_K=\frac{\sig_{n-2}}{n-1}\, M \rho_L^{n-1},\ee where $M$ is the
Minkowski-Funk transform (\ref{mf}). We denote by $\I\B_n^\infty$
the  subclass of $\I\B_n$ consisting of infinitely smooth bodies.
Intersection bodies of centered convex bodies first appeared in the
work of Busemann \cite{Bu} who did not gave them a particular name.
A more general class $\I_n$ of star bodies, which was called {\it
the class of intersection bodies} (without wording ``of star
bodies") was defined in \cite{GLW} as a collection of all bodies
$K\in\bbk^n$ with the property \be\label{eqvm} \rho_K =M \mu\quad
\text{for some} \quad \mu \in \M_+(S^{n-1}).\ee
  By Definition
 \ref{df}, the class $\I_n$ is a member
 of the family $\{\K_{\a,n}\}$ corresponding to $\a=1$. Thus, Theorems
 \ref{tcl} and \ref{bpd} imply the following known statement.
\begin{theorem}\label{bpd1} A body  $K \in\bbk^n$ is an intersection body, i.e., $K \in
\I_n$ if and only if $||\cdot||_K^{-1}$ is a positive definite
distribution. The class $\I_n$ is the closure of the classes
$\I\B_n$ and $\I\B_n^\infty$ of  intersection bodies of star bodies
in the radial metric.
\end{theorem}

The first part of this theorem can be found in \cite[Theorem
4.1]{K3}. Regarding the second part, see \cite[Theorem 5.5]{GZ},
 \cite{GLW}, and references therein.

{\bf Example 2.}\label{ex2} The following  extension of the
definitions in Example 1 to sections of arbitrary dimension $0<i<n$
was suggested by Koldobsky \cite{K2}. According to \cite{K2}, a body
$K\in \bbk^n$ is an $i$-intersection body of a body $L\in \bbk^n$
(we write $K=\I\B_i (L)$) if \be\label{ib3} \vol_{i} (K\cap
\xi)=\vol_{n-i} (L\cap \xi^\perp)\quad \forall \xi \in G_{n,i}, \ee
or, in other words, \be\label{ib4} \frac{\sig_{i-1} }{i}(R_{i}
\rho_K^{i})(\xi)=\frac{\sig_{n-i-1}}{n-i}(R_{n-i,\perp}
\rho_L^{n-i})(\xi)\quad \forall \xi \in G_{n,i}. \ee We denote by
$\I\B_{i,n}$ the class of all star bodies  with this property and by
$\I\B_{i,n}^\infty$ the subclass of $\I\B_{i,n}$ consisting of
infinitely smooth bodies.
 The above definition has a remarkable symmetry: \be\label{sym}
 K=\I\B_i (L)  \Longleftrightarrow L=\I\B_{n-i} (K).\ee
 We generalize this definition as follows.
\begin{definition}\label{def4}
 A body $K\in \bbk^n$ is  an
 $i$-intersection body if there is a measure $\mu\in \M_+(S^{n-1})$
  such that
\be\label{ib5} R_{i} \rho_K^{i}=R_{n-i,\perp} \mu. \ee We denote by
$\I_{i,n}$ the class of $i$-intersection bodies in $\bbr^n$.
\end{definition}
The  concept of the $i$-intersection body based on (\ref{ib3}) was
introduced by Koldobsky \cite[Definition 3]{K4}. His definition  is
given in terms of the Fourier transforms. Our Definition \ref{def4}
agrees with (\ref{eqvm}) and uses the language of Radon transforms.
As we shall see below, both definitions are equivalent.
\begin{theorem}\label{thm}\

{\rm(i)} For $\a=i$, the classes $\{\K_{\a,n}\}$,
 $\{\K\B_{\a,n}\}$, and    $\{\K\B_{\a,n}^\infty\}$
 coincide with $\{\I_{i,n}\}$, $\{\I\B_{i,n}\}$,
and $\{\I\B_{i,n}^\infty\}$, respectively.

{\rm(ii)} A body  $K \in\bbk^n$ is an $i$-intersection body, i.e.,
$K \in \I_{i,n} $, if and only if $||\cdot||_K^{-i}$ is a positive
definite distribution.

{\rm(iii)} The class $\I_{i,n}$ is the closure of the classes
$\{\I\B_{i,n}\}$ and $\{\I\B_{i,n}^\infty\}$ of $i$-intersection
bodies of star bodies in the radial metric.

{\rm(iv)} If a  body  $K \in\bbk^n$ is infinitely smooth and
$K=\I\B_i (L)$, then
 \be\label{kl}
\rho_L^{n-i}=\frac{\pi^{i-n/2}(n-i)}{i}\,M^{1-n+i} \rho_K^{i}.\ee
\end{theorem}
\begin{proof} (i) Let us prove that  equality (\ref{ib5}) is equivalent to
\be\label{cti}\rho_K^{i}=\tilde c^{-1}\,M^{1-i}\mu, \qquad \tilde c=
\frac{\pi^{i-n/2}\, \sig_{n-i-1}}{\sig_{i-1}}.\ee Indeed,
(\ref{ib5}) means that for any test function $\vp \in \D(G_{n,i})$,
we have $(\rho_K^{i}, R_i^*\vp)=(\mu, R^*_{n-i}\vp^\perp)$, or, by
(\ref{ccrl}), $(\rho_K^{i}, R_i^*\vp)=\tilde c^{-1}(\mu, M^{1-i}
R_i^* \vp)$. Since any function $\om \in \D_e(S^{n-1})$ can be
expressed as $\om=R_i^* \vp$ for some  $\vp \in \D(G_{n,i})$ (see
Lemma \ref{l34}), we are done, i.e., $\{\K_{i,n}\}=\{\I_{i,n}\}$.
The equalities $\{\K\B_{\a,n}\}=\{\I\B_{i,n}\}$ and
$\{\K\B_{\a,n}^\infty\}=\{\I\B_{i,n}^\infty\}$ can be proved
similarly: just use (\ref{ib4}) and replace $\mu$ by
$c\,\rho_L^{n-i}$ with $c=i\sig_{n-i-1}/(n-i)\sig_{i-1}$.

Statements (ii) and (iii) follow from Theorems   \ref{bpd} and
\ref{tcl}, respectively.

(iv) We make use of (\ref{conn}) with $\a=1-n+i$ and $f=\rho_K^{i}$.
Owing to (\ref{lim}), it can be written in the form \be \label {bu1}
R_i \rho_K^{i}=\tilde c R_{n-i,\perp}M^{1-n+i} \rho_K^{i}\ee with
constant $\tilde c$ as  in (\ref{cti}). On the other hand, if $K$ is
an infinitely smooth body and $ K =\I\B_i (L)$, then, by
(\ref{ib4}), \be \label {bu2} R_{i}
\rho_K^{i}=\frac{i\,\sig_{n-i-1}}{(n-i)\,\sig_{i-1} }\,R_{n-i,\perp}
\rho_L^{n-i}.\ee Comparing (\ref{bu1}) and (\ref{bu2}), we obtain
(\ref{kl}).
\end{proof}

Some comments are in order.

1. Statement (ii) of Theorem \ref{thm} was proved by Koldobsky (see
\cite{K2}, Theorem 4) for the case, when $K$ is an $i$-intersection
bodies of a star body with $C^\infty$-boundary. In our notation it
means that $K\in \I\B_{i,n}^\infty$. This result was extended in
\cite{K4} to a more general class of $i$-intersection bodies, which
coincides with $ \I_{i,n}$. Our argument essentially differs from
that in cited works.

2. Unlike the case $i=1$, when a body $ K=\I\B(L)$ can be
constructively realized
 for {\it any} origin-symmetric star body $L$, this is
not so if $i>1$, when the definition of $\I\B_i(L)$ is purely
analytic. It is known that for $i>3$, owing to Theorem \ref{thm} and
the symmetry (\ref{sym}),  the $i$--intersection body $\I\B_i (L)$
is not defined if $L$ is the unit ball $\{x \in \bbr^n \, : \,
|x_1|^4 + \dots |x_n|^4 <1\}$ of the space $\ell^n_4$. The reason is
that $||x||_L^{i-n}$ is not a positive definite distribution; see
Theorem 2 in \cite {K5}.

{\bf Problem.} The last remark provokes the following question: Is
 there a reasonable way to generalize the original Busemann's  construction
\cite{Bu} to
 sections of codimension $>1$ so that it would be
 applicable to {\it every} origin-symmetric star body? Clearly,  the normal
 unit vector in this construction should be replaced by the   orthonormal
 frame, the
 element of the corresponding Stiefel manifold.

{\bf Example 3.} For each origin-symmetric convex body $K$ in
$\bbr^n$, the Minkowski functional  $ ||\cdot||_K$ is a norm in
$\bbr^n$. A well known result going back to P. L\'evy \cite[Section
6.1]{K3} says that the space $(\bbr^n,||\cdot||_K)$ embeds
isometrically into $L_p, \; p>0$, if and only if there is a measure
$\mu \in \M_+ (S^{n-1})$ such that \be\label{lev}
||\theta||_K^p=\intl_{S^{n-1}} |\theta\cdot u|^p \, d\mu (u).\ee For
$p\neq 2,4, \ldots$, equality (\ref{lev}) is obviously equivalent to
$K \in \K_{\a,n}$ with $\a=-p$. Passing from $\a=-p$ to arbitrary
$\a$ in Definition \ref{df} allows us to extend the wording ``embeds
isometrically into $L_p$" to negative  $p$; see \cite{K3}. For
$p\neq -n, -n-2, \ldots$, this is equivalent to $K \in \K_{-p,n}$.

\subsection{Zhang's class of intersection bodies}

The following class of bodies was  introduced by Zhang \cite{Z2} in
his research related to the lower dimensional Busemann-Petty
problem.
\begin{definition}\label{def3} An origin-symmetric star body $K$ in $\bbr^n$ is called an
$i^*$-intersection body if there is a measure $\nu \in
\M_+(G_{n,i})$ such that  $\rho_K^{n-i}=R^*_i \nu$.
\end{definition}
Here abbreviation $i^*$  has been chosen  to distinguish this class
of bodies from that in Definition \ref{def4} and to indicate
implementation of the dual Radon transform $R^*_i$. Another notation
for both classes  of bodies was utilized in \cite{Mi}.

Connection between the class of $i^*$-intersection bodies and that
of $i$-intersection bodies is an important  problem intimately
related to the lower dimensional Busemann-Petty problem for sections
of convex bodies (the case of $2$ and $3$-dimensional sections is of
primary interest). To the best of our knowledge, both problems are
still open; see \cite{K3} for details.
 The fact that each  $(n-i)^*$-intersection body  is an
$i$-intersection body is due to Koldobsky \cite[Corollary 3]{K4},
who proved it using the Fourier transform technique and isometric
embedding in $L^p$-spaces. The theorem below includes this statement
and provides aditional information. Our proof differs from that in
\cite{K4} and does not involve the Fourier transform.

 Given a measure $ \nu \in \M_+(G_{n,n-i})$,  we
define a measure $ \nu^\perp \in \M_+(G_{n,i})$ by $ (\nu^\perp,
\vp)=(\nu, \vp^\perp)$ where $\vp \in C(G_{n,n-i}),\; \vp^\perp
(\xi)=\vp(\xi^\perp), \; \xi \in G_{n,i}$.

\begin{theorem}\label{th10} Each  $(n-i)^*$-intersection body $K$ is an
$i$-intersection body.  Specifically, if $\rho_K^{i}=R_{n-i}^*
\nu,\; \nu\in \M_+(G_{n,n-i})$, then $R_{i}
\rho_K^{i}=R_{n-i,\perp}\mu$ where $\mu= R_i^* \nu^\perp$.
Conversely, if $K$ is an $i$-intersection body, i.e., $R_{i}
\rho_K^{i}=R_{n-i,\perp}\mu$ for some $\mu\in \M_+(S^{n-1})$, then
$K$ is an  $(n-i)^*$-intersection body provided that $\mu $ is
represented in the form $\mu= R_i^* \nu^\perp$ for some $\nu\in
\M_+(G_{n,n-i})$. If the latter is true, then $\rho_K^{i}=R_{n-i}^*
\nu$.
\end{theorem}
\begin{proof}  Let $K$ be an $(n-i)^*$-intersection body, i.e.,
$\rho_K^{i}=R_{n-i}^* \nu,\; \nu\in \M_+(G_{n,n-i})$. For any
function $\psi \in \D(G_{n,i})$,  we have \bea (R_{i} \rho_K^{i},
\psi)&=& (\rho_K^{i},R_{i}^*\psi)=(R_{n-i}^* \nu,R_{i}^*\psi)=(\nu,
R_{n-i}R_{i}^*\psi)\nonumber \\
&=&(\nu^\perp,R_{n-i, \perp}R_{i}^*\psi) \qquad \text{\rm (use (\ref{kr}))}\nonumber \\
&=&\tilde c^{-1}(\nu^\perp,R_i M^{1-i}R_{i}^*\psi).\nonumber\eea If
we set $\mu= R_i^* \nu^\perp$, then the last expression becomes
$\tilde c^{-1}(\mu, M^{1-i}R_{i}^*\psi)$ which coincides with
$(\mu,R_{n-i}^*\psi^\perp)= (R_{n-i,\perp}\mu, \psi)$ by
(\ref{ccrl}). This gives the result.

Conversely, let $R_{i} \rho_K^{i}=R_{n-i,\perp}\mu$, $\mu\in
\M_+(S^{n-1})$. We take a test function $\om\in \D_e (S^{n-1})$ and
represent it in the form $\om=R_{i}^*\psi$, $\psi \in \D(G_{n,i})$
(this is possible by Lemma \ref{l34}). Then $$ (\rho_K^{i},
\om)=(\rho_K^{i}, R_{i}^*\psi)=(R_{i} \rho_K^{i},
\psi)=(R_{n-i,\perp}\mu,\psi).$$ Since $\mu= R_i^* \nu^\perp$, we
 continue \bea&=&(R_i^* \nu^\perp,R_{n-i}^*\psi^\perp)=(\nu,
R_{i,\perp}R_{n-i}^*\psi^\perp)\nonumber \\&&\text{\rm (use (\ref{kr}) with $i$ replaced by $n-i$)}\nonumber \\
&=&\tilde c^{-1}(\nu,R_{n-i} M^{1-n+i}R_{n-i}^*\psi^\perp)=\tilde
c^{-1}(R_{n-i}^*\nu,M^{1-n+i}R_{n-i}^*\psi^\perp).\nonumber \eea By
(\ref{ccrl}), this coincides with
$(R_{n-i}^*\nu,R_i^*\psi)=(R_{n-i}^*\nu,\om)$. Hence
$\rho_K^{i}=R_{n-i}^* \nu$, and the proof is complete.
\end{proof}

\end{document}